\documentclass[10pt]{article}

\usepackage{theorem,amssymb,amsmath,latexsym}

\usepackage{graphicx}

\usepackage{color}                    
\definecolor{red}{rgb}{1,0,.2}        
\newcommand*{\clrred}[1]{{\color{red} #1}} 
\newcommand{\cred}[1]{\clrred{ #1}}   
\definecolor{cjp}{rgb}{.1,.7,.2}        
\definecolor{fmdc}{rgb}{1,0,.8}        

\newcommand{\TotKap}{{\rm T\!ot^\kappa}}
\newcommand{\TotNu}{{\rm T\!ot^\nu}}
\newcommand{\TotFIB}{{\rm T\!ot^{\rm \sc FIB}}}
\newcommand{\redcdot}{\cred{\cdot}}

\newtheorem{theorem}{Theorem}[section]

\newtheorem{lemma}[theorem]{Lemma}

\newtheorem{proposition}[theorem]{Proposition}
\theorembodyfont{\rm}

\newtheorem{remark}[theorem]{Remark}
\newtheorem{example}[theorem]{Example}

\advance \topmargin by -\headheight
\advance \topmargin by -\headsep
\evensidemargin \oddsidemargin
\begin{document}

\begin{centering}
{\Large \textbf{Counting base phi representations}}

\bigskip

{\bf \large Michel Dekking and Ad van Loon}

\bigskip

{\footnotesize  \it Adresses M.~Dekking:  CWI, Amsterdam and Delft University of Technology, Faculty EEMCS.}

\medskip

{\footnotesize \it Email adresses:  Michel.Dekking@cwi.nl, f.m.dekking@tudelft.nl, advloon@upcmail.nl}

\bigskip

{ \bf  April 22, 2023}

\end{centering}

\medskip

\begin{abstract}
 \noindent In a base phi representation a natural number is written as a sum of powers of the golden mean $\varphi$.   There are many ways to do this. How many?
 Even if the number of powers of $\varphi$ is finite, then any number has infinitely many base phi representations. By not allowing an expansion to end with the digits 0,1,1, the number of expansions becomes finite, a solution proposed by Ron Knott. Our first result is a recursion to compute this number of expansions. This recursion is closely related to the recursion given by Neville Robbins to compute the number of Fibonacci representations of a number, also known as Fibonacci partitions. We propose another way to obtain finitely many expansions, which we call the natural base phi expansions. We prove that these are closely connected to the  Fibonacci partitions.
\end{abstract}

\medskip


\bigskip


\date{\today}

\bigskip

\section{Introduction}

A natural number $N$ is written in base phi if $N$ has the form
  $$N= \sum_{i=-\infty}^{\infty} a_i \varphi^i,\vspace*{-.0cm}$$
where  $a_i=0$ or $a_i=1$, and where $\varphi:=(1+\sqrt{5})/2$ is the golden mean.

There are infinitely many ways to do this.
When the number of powers of $\varphi$ is finite we  write these representations (also called expansions) as
  $$\alpha(N) = a_{L}a_{L-1}\dots a_1a_0\redcdot a_{-1}a_{-2} \dots a_{R+1}a_R.$$

\medskip

 Infinitely many expansions can be generated in a rather trivial way from expansions with just a few powers of $\varphi$ using the replacement $100\rightarrow011$ at the end of the expansion.

 So we use Knott's truncation rule from \cite{Knott2}:
 \begin{equation}\label{eq:Knott}
 a_{R+2}a_{R+1}a_{R}\ne 011.
  \end{equation}
 Let $\TotKap(N)$ be the  number of base phi expansions of the number $N$ satisfying Equation (\ref{eq:Knott}):

\medskip

$\TotKap =\:0, 1, 1, 2, 3, 3, 5, 5, 5, 8, 8, 8, 5, 10, 13, 12, 12, 13, 10, 7, 15, 18, 21, 16, 20, 20, 16, 21, 18, 15, 7, 17,\ldots$\footnote{In OEIS (\cite{oeis}):  A289749 \;		Number of ways not ending in 011 to write n in base phi.}	

\medskip

In 1957 George Bergman (\cite{Bergman}) proposed restrictions on the digits $a_i$ which entail that the representation becomes unique and finite.
This is generally accepted as \emph{the} representation of the natural numbers in base phi.
A natural number $N$ is written in the Bergman representation  if $N$ has the form
  $$N= \sum_{i=-\infty}^{\infty} d_i \varphi^i,\vspace*{-.0cm}$$
  with digits $d_i=0$ or $d_i=1$, and where $d_{i+1}d_i = 11$ is not allowed.  We  write these representations as
  $$\beta(N) = d_{L}d_{L-1}\dots d_1d_0\redcdot d_{-1}d_{-2} \dots d_{R+1}d_R.$$

  \bigskip

  A natural number $N$ is written in base Fibonacci if $N$ has the form  
  $$N= \sum_{i=2}^{\infty} c_i F_i,\vspace*{-.0cm}$$
 where $c_i=0$ or  $c_i=1$, and $(F_i)_{i\ge 0}=0,1,1,2,3,\ldots$ are the Fibonacci  numbers.

 \smallskip

 Let $\TotFIB(N)$ be the  total number of Fibonacci expansions of the number $N$. Then

$$\TotFIB=\:1, 1, 1, 2, 1, 2, 2, 1, 3, 2, 2, 3, 1, 3, 3, 2, 4, 2, 3, 3, 1, 4, 3, 3, 5,\ldots\footnote{In OEIS (\cite{oeis}): A000119\;Number of representations of n as a sum of distinct Fibonacci numbers.}$$
This sequence has received a lot of attention, see e.g., the papers \cite{Klarner66}, \cite{Klarner68}, \cite{Carlitz-1968}, \cite{car-sco-hog}, \cite{Robbins}, \cite{Bicknell}, \cite{Stockmeyer}, and \cite{Chow}.

In 1952 the paper \cite{Lekker} proposed restrictions on the digits $c_i$ which entail that the representation becomes unique.
This is known as the Zeckendorf expansion of the natural numbers after the paper \cite{Zeck}.

A natural number $N$ is written in the Zeckendorf representation  if $N$ has the form
  $$N= \sum_{i=2}^{\infty} e_i F_i,\vspace*{-.0cm}$$
  with digits $e_i=0$ or $e_i=1$, and where $e_{i+1}e_i = 11$ is not allowed.

The Fibonacci representation and the base phi representation are  closely related. We make a list.

\bigskip

\begin{tabular}{|c|c|c|}
  \hline
  Property & Fibonacci & Base phi  \\
  \hline\\[-.4cm] \hline
   & $F_n:$ \quad $n\ge 2$& $\varphi^{n}:\quad n$ integer      $\phantom{\Big (}$ \\
  \hline
  Fundamental recursion & $F_{n+1}=F_n+F_{n-1}$ & $\varphi^{n+1}=\varphi^n+\varphi^{n-1} $  $\phantom{\bigg (}$  \\
  Golden mean flip & $100\rightarrow 011$ & $100\rightarrow 011$     \\
  \hline
  Unique expansion & Zeckendorf & Bergman        \\
  Condition on the digits & no 11 & no 11         \\
  \hline
  Fundamental intervals & $[F_n, F_{n+1}-1]$ & $[L_{2n},\,L_{2n+1}],\,[L_{2n+1}+1,\, L_{2n+2}-1]$      \\
  Examples $F_5=5$, $L_4=7$& $[5,7]= [2\,2\,1]$ & $ [7,11] = [5\,8\,8\,8\,5] $   \\
  Examples $F_6=8$, $L_5=11$ & $[8,12]=[3\,2\,2\,3\,1]$ & $[12,17]=[10\, 13\, 12\, 12\, 13\, 10]$   \\
  \hline
\end{tabular}

\bigskip

 Here the $L_n$ are the Lucas numbers defined by $L_0=2, L_1=1$ and $L_{n+1}=L_n+L_{n-1}$ for $n\ge1$.\\
  The intervals $\Lambda_{2n}=[L_{2n},\,L_{2n+1}]$, $\Lambda_{2n+1}=[L_{2n+1}+1,\, L_{2n+2}-1]$ are called the even and odd Lucas intervals.

\medskip

Replacing the digits 100 in an expansion by 011 will be called a \emph{golden mean flip}. Our Theorem \ref{thm:fromB} shows that any base phi expansion can be obtained from the Bergman expansion by a finite number of such golden mean flips. There is a special case which needs attention, which we illustrate with an example. Let $N=4$. Then $\beta(4)=101\redcdot 01$.
Applying the golden mean flip at the right gives the expansion $101\redcdot0011$, which is not an allowed expansion.
However, if we apply a second golden mean flip we can obtain    $100\redcdot 1111$, which \emph{is} an allowed expansion. We call this operation a \emph{double golden mean flip}.

\medskip

In Section \ref{sec:KnotTot} we determine a formula for $\TotKap(N)$. In Section \ref{FibLuc} we give simple formula's for $N=F_n$, and for $N=L_n$. In Section \ref{sec:TotNu} we introduce a new way to count expansions, by defining \emph{natural expansions}, and give a formula for $\TotNu(N)$,  the number of natural base phi expansions of $N$. We moreover show that $(\TotNu(N))$ is a subsequence of the sequence of total numbers of Fibonacci representations. Section \ref{sec:compare} gives important information on the different behaviour of phi expansions on the odd and the even Lucas intervals.

\section{A recursive formula for the number of Knott expansions}\label{sec:KnotTot}

In this section we determine a formula for $\TotKap(N)$ for each natural number $N$.

\medskip

Because of the fundamental recursion $\varphi^{n+1}=  \varphi^{n}+\varphi^{n-1}$, one can transform a
base phi expansion $\alpha(N) = a_{L}a_{L-1}\dots a_1a_0\redcdot a_{-1}a_{-2} \dots a_{R+1}a_R$ of $N$ with $a_{i+1}a_{i}a_{i-1}=100$ to another base phi expansion of $N$, by the map
$$ T_i: a_{i+1}a_{i}a_{i-1} \rightarrow [a_{i+1}-1][a_{i}+1][a_{i-1}+1] ,$$
where $R-1\le i\le L-1$. This is a more detailed definition of the golden mean flip.\\In the definition we put of course $a_{R-1}=a_{R-2}=0$.

The map $T_i$ has an inverse denoted $U_i$ for $R-1\le i\le L$ given by
$$ U_i: a_{i+1}a_{i}a_{i-1}   \rightarrow [a_{i+1}+1][a_{i}-1][a_{i-1}-1] ,$$
as soon as $a_{i+1}a_{i}a_{i-1}=011$.  We call this map the  {\it reverse golden mean flip}.

Example: $U_{1}(110\redcdot 01)=1000\redcdot01$.

\begin{theorem}\label{thm:fromB}
  Any finite base phi expansion $\alpha(N)$ with digits 0 and 1 of a natural number $N$ can be obtained from the
     Bergman expansion $\beta(N)$ of $N$ by a finite number of applications of the  golden mean flip.
\end{theorem}

\noindent  {\it Proof:} We prove this by showing that any expansion of $N$ will be mapped to its Bergman expansion by a finite number of applications of the reverse golden mean flip. Let $\alpha(N) = a_{L}a_{L-1}\dots a_1a_0\redcdot a_{-1}a_{-2} \dots a_{R+1}a_R$ be an expansion of $N$ with digits 0 and 1.
When 11 does not occur in $\alpha(N)$, then $\alpha(N)=\beta(N)$, and there is nothing to do. Otherwise, let $m:=\max\{i: a_{i}a_{i-1}=11\}$. First, suppose $m\le L-2$. Then by the definition of $m$, we have $a_{i+1}=0$. So for the two possibilities $a_{i+2}=0$ and $a_{i+2}=1$
$$ U_{i}(\ldots 0a_{i+1}a_{i}a_{i-1}\ldots)=U_{i}(0011)=0100, \quad {\rm and}\; U_{i}(\ldots 1a_{i+1}a_{i}a_{i-1}\ldots)=U_{i}(1011)=1100.$$
Note that in the first case the total number of 11 occurring in the expansion of $N$ has decreased by 1, and in the second case it remained constant.
However, in the second case the $m$ of $U_{i}(\alpha(N))$ has increased by 2. If we keep iterating the   reverse golden mean flip on the left most occurrence of 11, then either 0011 will occur, or if not, then $\alpha(N)=1101\dots$. This is the case $m=L$, where there \emph{is} a decrease in the number of 11, since $U_{L}(1101\dots)=10001\dots$.
Conclusion: in all cases the number of 11 will decrease by at least 1 after a finite number of applications of the  reverse golden mean flip. So after a finite number of applications of the   reverse golden mean flip we reach an expansion with no occurrences of 11. By definition, this is the Bergman expansion.

The case $m=L$ has already been considered above, the case $m=L-1$ corresponds to $\alpha(N)=011\dots$, where an application of the   reverse golden mean flip leads also to a decrease in the number of 11. \hfill $\Box$

\begin{example} Let $N=5$, with $\beta(5)=1000\redcdot1001$:
$$10\underline{1}\redcdot1111\rightarrow \underline{1}10\redcdot0111\rightarrow 1000\redcdot0\underline{1}11\rightarrow 1000\redcdot1001, \quad {\rm with\: maps\;} U_0, U_{2}, U_{-2}.$$
\end{example}

 Our proof for $\TotKap$ resembles  the work of Neville Robbins \cite{Robbins} on Fibonacci representations, but we have to incorporate the double golden mean flip defined in the Introduction. It then appears that the two recursions for Fibonacci representations and golden mean (Knott) representations are the same, but that there is a difference in the initial conditions.

\medskip

The emphasis will be on the manipulation of $0$-$1$-words, not on numbers.\\
Let $\beta(N) = d_{L}d_{L-1}\dots d_1d_0\redcdot d_{-1}d_{-2} \dots d_{R+1}d_R$. By removing the radix point, we obtain a $0$-$1$-word $B(N):=d_{L}d_{L-1}\dots d_1d_0d_{-1}d_{-2} \dots d_{R+1}d_R.$
Let us denote  $r(B(N)):=\TotKap(N)$.

More generally, $r(w)$ is the number of words satisfying the Knott condition that can be obtained from a word $w$ by golden mean flips. Note that  in general the representations that we obtain are not representations of a natural number---not for any choice of the radix point. An example is $w=100001$, which represents $\varphi^5+1$. Nevertheless, these words represent numbers $a+b\varphi$ with non-negative natural numbers $a$ an $b$ in the ring $\mathbb{Z}(\varphi)$.

 For example $w=100001$ represents $5\varphi+4$. This is the justification for  continuing with the terminology of representations.

\medskip

Here are two basic examples.\\[-.7cm]
\begin{align}\label{eq:basic}
 r(10^s) = & \:\frac12 s +1 \quad   s\;{\rm even} \\
 r(10^s) = & \:\frac12(s+1) \quad   s\;{\rm odd}
\end{align}
This follows easily by making golden mean flips from left to right.

\medskip
Suppose the Bergman representation $\beta(N)$ of a number $N$ contains $n+1$ ones. Then we can write for some numbers $s_1,s_2,\ldots,s_n$

$$B(N) = 10^{s_n}\ldots 10^{s_2}\, 10^{s_1}\,1.$$
We start with the case $n=2$, so $$B(N) = 10^{s_2}10^{s_1}\,1.$$
Let us call $I_2:=10^{s_2}$ the \emph{initial segment} of $B(N)$, and $T_1:=10^{s_1}\,1$ the \emph{terminal segment} of $B(N)$.\\
We want to deduce $r(B(N))=r(I_2T_1)$ from the number of representations $r(I_2)$ and $r(T_1)$.
There are two cases to consider.

\medskip

\noindent Type 1:\; Arbitrary combinations of representations of $I_2$ and $T_1$.\\
Type 2:\; Arbitrary combinations of representations of $I_2$ and $T_1$ \emph{plus} an `overlap' combination.

\medskip

Type 1 typically occurs if $s_2$ is even. For example for the case $s_2=4$, we have the three representations $10000,\:  01100,\: 01011$. Note that in general these representations always end in 00 or 11.

So for Type 1 one has simply
\begin{equation}
r(B(N))=r(I_2T_1)=r(I_2)r(T_1).
\end{equation}

But for $s_2$ odd, for example $s_2=5$,    $100000,\: 011000,\:  010110$ are the three representations of $I_2$. Note that in general these representations always end in 00 or 10.\\
So if a representation $0v$ of $T_1$ starts with a 0, then the representation $ 010110\,0v$ generates an  `overlap' representation $ 010101\,1v$ via the golden mean flip.\\
Obviously it is true in general that an $I_2$ word with $s_2$ odd will have exactly one representation that ends in 10. Also important: there is no representation that ends in 01.
Therefore, if $r^{(i)}(T_1)$ denotes  the number of representations of $T_1$ starting with $i$ for $i=0,1$, then we obtain for Type 2:
\begin{equation}\label{eq:rBN1}
r(B(N))=r(I_2T_1) =r(I_2)r(T_1)+r^{(0)}(T_1).
\end{equation}

It thus follows from Equation (\ref{eq:rBN1}), the trivial equation  $r^{(0)}(T_1) +  r^{(1)}(T_1)=r(T_1)$, and the fact that the segment $T_1=10^{s_1}\,1$ has just one representation that starts with a 1, that
\begin{equation}\label{eq:rBN2}
r(B(N))=r(T_1)[r(I_2)+1]-r^{(1)}(T_1)= r(T_1)[r(I_2)+1]-1.
\end{equation}

\medskip

We continue with the case $n=3$, so $$B(N) = 10^{s_3}10^{s_2}10^{s_1}\,1.$$
Now $I_3:=10^{s_3}$ is the \emph{initial segment}, and $T_2:=10^{s_2}10^{s_1}\,1$ the \emph{terminal segment}.\\
As before there are two cases to consider  to compute $r(B(N))=r(I_3T_2)$.

\medskip

\noindent Type 1:\; Arbitrary combinations of representations of $I_3$ and $T_2$.\\
Type 2:\; Arbitrary combinations of representations of $I_3$ and $T_2$ \emph{plus} an `overlap' combination.

\medskip

For Type 1 one has simply
\begin{equation}
r(B(N))=r(I_3T_2)=r(I_3)r(T_2).
\end{equation}

For Type 2 one has :
\begin{equation}\label{eq:rBN4}
r(B(N))=r(I_3T_2) = r(I_3)r(T_2)+r^{(0)}(T_2).
\end{equation}

Next, we split $T_2=I_2T_1$, where $I_2=:10^{s_2}$. Then we have, since $I_2$ has just one representation that starts with a 1, that $r^{(1)}(T_2)=r(T_1)$.
It thus follows from Equation(\ref{eq:rBN4}) and  $r^{(0)}(T_2) + r^{(1)}(T_2)= r(T_2)$ that      
\begin{equation}
r(B(N))=r(I_3)r(T_2)+r(T_2)-r^{(1)}(T_2)=r(T_2)[r(I_3)+1]-r^{(1)}(T_2)= r(T_2)[r(I_3)+1]-r(T_1).
\end{equation}

In the same way one proves that for any $k= 1,\ldots n-1$ the following formula holds for $s_{k+1}$ odd, where $T_{k+1}$ is split as $T_{k+1}=I_{k+1}T_k$.

\begin{equation}
r(T_{k+1})= r(T_k)[r(I_{k+1})+1]-r(T_{k-1}).
\end{equation}

Defining $r_{n}:=r(B(N)),\, r_k:=r(T_{k})$ for $k=1,\ldots,n-1$ and $r_0=1$ (cf. Equation (\ref{eq:rBN2})), we have obtained a recursion that computes $r(B(N))$.

\begin{theorem}\label{th:totalkappa} For a natural number $N$ let the Bergman expansion of $N$ have $n+1$ digits $1$. Suppose $\beta(N)= 10^{s_n}\ldots 10^{s_1}\,1$.
Let $\TotKap(N) = r_n$ be the number of Knott representations of  $N$.
Define the initial conditions: $r_0=1$ and $r_1 = \frac12 s_1 + 1$ if $s_1$ is even, $r_1=  \frac12(s_1\!+\!1)+1$ if $s_1$ is odd. Then for $n\ge 2$:
\begin{displaymath}
r_{n} = \begin{cases}
[\frac12s_n\!+\!1]\,r_{n-1}\hspace*{2.7cm}{\rm if\;} s_n \; {\rm is\: even}\\[.1cm]
 [\frac12(s_n\!+\!1)+1]r_{n-1}-r_{n-2} \qquad {\rm if\;} s_n \; {\rm is\: odd}\\[-.1cm]
\end{cases}
\end{displaymath}
\end{theorem}

\medskip

The initial condition for $r_1$ is different from the Fibonacci case: if $s_1$ is odd, then the base phi expansion has an extra representation that is generated by the ``double golden mean flip"  (see Section \ref{sec:compare}).

\section{Expansions of the Fibonacci numbers and  the Lucas numbers}\label{FibLuc}

Let $(F_n)=0,1,1,2,3,5,\dots$ be the Fibonacci numbers. We will determine the number of Knott representations of these numbers.
Then we first have to find a formula for the Bergman expansions of the Fibonacci numbers.
Recall that $L(N)$ is the left most position of a 1 in $\beta(N)$, and that $B(N)$ is $\beta(N)$ without the radix point in the expansion.
In the following proposition and its proof the simple 1-to-1 correspondence between the pair $(L(N), B(N))$ and the Bergman expansion $\beta(N)$ plays an essential role.

\begin{proposition}\label{prop:Fib}
a)  For $n\ge 3$ one has $L(F_n)= n-1$.\\
b) If $n\ge 3$ is odd, then $B(F_n)=(1000)^p1001$, with $p=(n-3)/2$.\\
c) If $n\ge 4$ is even, then  $B(F_n)=(1000)^p10001$, with $p=(n-4)/2$.
\end{proposition}

\noindent  {\it Proof:}
This will, of course, be proved by induction.
It is simple to check that $\beta(F_3)=\beta(2)=10\redcdot01, \beta(F_4)=\beta(3)=100\redcdot01, \beta(F_5)=\beta(5)=1000\redcdot1001$.
So the statements hold for $n=3,4,5$. We start the induction at $n=6$. Since $F_6=F_4+F_5$, we have\\[-.8cm]
\begin{align*}
  \beta(F_4)=  &\;\;\;\, 100\redcdot01 \\
 \beta(F_5)=  & \;\;1000\redcdot1001 \\
\beta(F_4)+\beta(F_5)=  & \;\; 1100\redcdot1101\\
\beta(F_4)+\beta(F_5)=  & \; 10001\redcdot0001.
\end{align*}
Here we applied the reverse golden mean flip twice in the last step. Since the last expansion does not have any 11, we must have $\beta(F_6)=10001\redcdot0001$, and $L(F_6)=5$.
Next we show what happens at $n=7$.
\begin{align*}
\beta(F_5)=  & \;\;\;\;1000\redcdot1001 \\
\beta(F_6)=  &\;\;\;10001 \redcdot0001 \\
\beta(F_5)+\beta(F_6)=  & \;\;\; 11001\redcdot1002\\
\beta(F_5)+\beta(F_6)=  & \; 100010\redcdot001001.
\end{align*}
Here we applied the reverse golden mean flip twice, and a shifted version of $\beta(2)=10\redcdot 01$ in the last step.
Since the last expansion does not have any 11, we must have $\beta(F_7)=100010\redcdot 001001$, and $L(F_7)=6$.

These addition schemes clearly generalize to $\beta(F_{n-2})+\beta(F_{n-1})$ with $n-2$ even, respectively odd, finishing the induction proof. \hfill $\Box$

\medskip

\begin{theorem}\label{th:TotFn}
For all $n\ge 1$ one has $\TotKap(F_n)=F_n$.
\end{theorem}

\noindent  {\it Proof:}  It is easily checked that the proposition holds for $n=1$ and $n=2$. So let $n\ge 3$.
According to Proposition \ref{prop:Fib}, the number of ones in $\beta(F_n)$ is $p+2$, with $p+2=(n+1)/2$ if $n$ is odd, and $p+2=n/2$ if $n$ is even.
Also, $\beta(F_n)=10^{s_{p+1}}\ldots 10^{s_k}\ldots 10^{s_1}\,1$, with $s_k=3$ for $k=2,\ldots, p+1$, and $s_1=2$ for $n$ odd, $s_1=3$ for $n$ even.

We apply Theorem \ref{th:totalkappa}. This yields that  $\TotKap(F_n)=r_{p+1}$,  the number of Knott representations of the Bergman representation of $F_n$
satisfies $$ r_{p+1}=3r_p-r_{p-1}. $$
Here the initial conditions are $r_0=1$, $r_1=s_1/2+1=2$ for $n$ even, and $r_1=(s_1+1)/2+1=3$ for $n$  odd.\\
Amusingly,  the same  recurrence relation holds for the subsequences of even and odd Fibonacci numbers:
\begin{equation}\label{eq:Fibrec}
F_{n+1}=F_{n}+F_{n-1}=2F_{n-1}+F_{n-2}=3F_{n-1}-F_{n-1}+F_{n-2}=3F_{n-1}-F_{n-3}.
\end{equation}

{\bf \rm (I)} \;Suppose  $n=2m+1$ is odd. Then  $p=m-1$, so $\TotKap(F_{2m+1})=r_m$.\\ We claim that $r_m=F_{2m+1}$ for all $m\ge 0$.\\
For $m=0$, we have $r_0=1 = F_1$, and for $m=1$ we have $r_1=2=F_3.$\\
For $m\ge 2$,
  $$ r_m = 3r_{m-1}-r_{m-2}=3F_{2m-1}-F_{2m-3}=F_{2m+1},$$
by the induction hypothesis and Equation (\ref{eq:Fibrec}).

{\bf \rm (II)} \; Suppose  $n=2m+2$ is even. Then  $p=m-1$, so $\TotKap(F_{2m+2})=r_m$.\\ We claim that $r_m=F_{2m+2}$ for all $m\ge 0$.\\
For $m=0$, we have $r_0=1 = F_2$, and for $m=1$ we have $r_1=3=F_4.$\\
For $m\ge 2$,
  $$ r_m = 3r_{m-1}-r_{m-2}=3F_{2m}-F_{2m-2}=F_{2m+2},$$
by the induction hypothesis and Equation (\ref{eq:Fibrec}).\\
Combining {\bf \rm (I)} and {\bf \rm (II)} yields the conclusion: $\TotKap(F_n)=F_n$ for all $n\ge 1$. \hfill $\Box$

\bigskip

At the Fibonacci numbers the total number of expansions is very large, but here we show that it is very small at the Lucas numbers $(L_n)$.

\begin{theorem}
  For all $n\ge 1$ one has $\TotKap(L_{2n})=\TotKap(L_{2n+1})=2n+1 $.
\end{theorem}

\noindent  {\it Proof:} The Lucas numbers have simple representations: $\beta(L_{2n}) = 10^{2n}\redcdot0^{2n-1}1,\; \beta(L_{2n+1}) = 1(01)^n\redcdot(01)^n$. This can be shown using  the  golden mean flip as in Theorem \ref{thm:fromB}.

So the representation  of $L_{2n}$ has only two ones. It follows therefore from Theorem \ref{th:totalkappa} that $\TotKap(L_{2n})=r_1=(s_1+1)/2+1=2n+1$, since $s_1=4n-1$ is odd.

The representation of $L_{2n+1}$ has $2n+1$ ones, and each $s_k$ of the blocks $10^{s_k}$ is equal to 1, which is odd. It follows therefore from Theorem \ref{th:totalkappa} that\; $\TotKap(L_{2n+1})=r_n=2r_{n-1}-r_{n-2}$.
And indeed, induction gives that $r_n=2(2n-1)-(2n-3)=2n+1$. \hfill $\Box$

\section{Natural base phi expansions}\label{sec:TotNu}

A consequence  of the application of the double golden mean shift is that length of the negative part of the Knott expansions may take two different values.

 To obtain what we will call the \emph{natural} expansions, let us delete all expansions that have a length of the negative part that is not equal to the length of the negative part of the Bergman expansion.

 For example in the case $N=4$ Knott proposes the three expansions $101\redcdot 01, 100\redcdot 1111$ and $11.1111 $. However, there is only one natural expansion: the Bergman expansion $101\redcdot 01$.

\medskip

Let $\TotNu(N)$ denote the number of natural base phi expansions. Then we have the following
\begin{align*}
  (\TotNu(N)) & = 1,1,2,2,1,5,5,4,5,4,3,1,10,13,12,12,13,10,6,11,12,\dots \\
    {\rm instead\: of} & \\
  (\TotKap(N)) & = 1,1,2,3, 3, 5, 5, 5, 8, 8, 8, 5, 10, 13, 12, 12, 13, 10, 7, 15, 18, \dots
\end{align*}

The number of natural base phi expansions can be determined in a way that is very similar to the Knott expansion case.

\begin{theorem}\label{th:totalnu} For a natural number $N$ let the Bergman expansion of $N$ have $n+1$ digits $1$. Suppose $\beta(N)= 10^{s_n}\ldots 10^{s_1}\,1$.
Let $\TotNu(N) = r_n$ be the number of natural base phi representations of  $N$.
Define the initial conditions: $r_0=1$ and $r_1 = \frac12 s_1+1$ if $s_1$ is even, $r_1=  \frac12(s_1\!+\!1)$ if $s_1$ is odd. Then for $n\ge 2$:
\begin{displaymath}
r_{n} = \begin{cases}
[\frac12s_n\!+\!1]\,r_{n-1}\hspace*{2.7cm}{\rm if\;} s_n \; {\rm is\: even}\\[.1cm]
 [\frac12(s_n\!+\!1)+1]r_{n-1}-r_{n-2} \qquad {\rm if\;} s_n \; {\rm is\: odd}\\[-.1cm]
\end{cases}
\end{displaymath}
\end{theorem}

\noindent  {\it Proof:} This follows directly from Theorem \ref{th:totalkappa} and its proof. The only difference between the process of generating all Knott expansions and all natural expansions is the double golden mean flip, which is performed in the Knott expansion at the  segment $10^{s_1}1$, and  only when $s_1$ is odd.
So  $\TotNu(N) = r_n$ satisfies the same recursion as $\TotFIB(N)$, except that  $r_1=  \frac12(s_1\!+\!1)+1$ has to be replaced by  $r_1=  \frac12(s_1\!+\!1)$ in the case that $s_1$ is odd.\hfill $\Box$

\medskip

We will determine the total number of natural expansions of the Fibonacci numbers. First we present a lemma that emphasizes the inter-connection between the Fibonacci and the Lucas numbers.
Recall the even and odd Lucas intervals  $\Lambda_{2n}=[L_{2n},\,L_{2n+1}]$, $\Lambda_{2n+1}=[L_{2n+1}+1,\, L_{2n+2}-1]$  (cf.~\cite{Dekk-phi-FQ}).

 \begin{lemma}\label{lem:FibLuc}
 For all $n=1,2,\ldots$ one has $F_{2n+2}\in \Lambda_{2n},\: F_{2n+3}\in \Lambda_{2n+1}$.
 \end{lemma}

\noindent  {\it Proof:} By induction. For $n=1$ we have $F_3=4\in \Lambda_2=[3,4]$, and   $F_5=4\in \Lambda_3=[5,6]$.

For $n=2$ we have $F_6=8\in \Lambda_4=[7,11]$, and $F_7=13\in \Lambda_5=[12,17]$.

\medskip

Suppose the statement of the lemma has been proved for $F_{2n+1}$ and $F_{2n+2}$. So we know
\begin{align*}
F_{2n+1} &\in [L_{2n-1}+1,\,L_{2n}-1]=\Lambda_{2n-1}\\
F_{2n+2} &\in [L_{2n},\, L_{2n+1}]\hspace*{1.2cm}=\Lambda_{2n.}
\end{align*}
Adding the numbers in these two equations vertically, we obtain
$$F_{2n+3} \in [L_{2n+1}+1,\, L_{2n+2}-1]=\Lambda_{2n+1}.$$
 We can then write\\[-.9cm]
 \begin{align*}
F_{2n+2} &\in [L_{2n},\,L_{2n+1}]\hspace*{1.55cm}= \Lambda_{2n}\\
F_{2n+3} &\in [L_{2n+1}+1,\, L_{2n+2}-1]=\Lambda_{2n+1.}
\end{align*}
This time, adding gives\\[-.9cm]
$$F_{2n+4} \in [L_{2n+2}+1,\, L_{2n+3}-1].$$
Since $F_{2n+4} \ne L_{2n+2}$, this implies that $F_{2n+4} \in \Lambda_{2n+2}$. \hfill $\Box$

\begin{theorem}\label{th:TotNuFn}  For all $n=0,1,2,\ldots$ one has $\TotNu(F_{2n+2})=F_{2n+1}$  and  $\TotNu(F_{2n+3})=F_{2n+3}$.
\end{theorem}

\noindent  {\it Proof:} We use the result from Proposition \ref{prop:s}, which gives that for all $N$ from $\Lambda_{2n+1}$  if $\beta(N)=...10^{s_1}1$, then $s_1$ is even. So for all $N$ from $\Lambda_{2n+1}$ we have that the total number of natural expansions is equal to the total number of Knott expansions. In particular we obtain from Lemma \ref{lem:FibLuc}, using Theorem \ref{th:TotFn}, that
$$\TotNu(F_{2n+3})=\TotKap(F_{2n+3})=F_{2n+3}.$$
From Proposition \ref{prop:Fib}, part c) we have that $B(F_{2n+2})=(1000)^p10001$ with $p=(2n+2-4)/2=n-1$. Therefore $r_n$ satisfies the recurrence relation $r_n=3r_{n-1}-r_{n-2}$, with $r_1=\frac12(3+1)=2=F_3$.
This is the recurrence relation for the Fibonacci numbers with odd indices, cf.~Equation (\ref{eq:Fibrec}). Therefore $r_n=F_{2n+1}$. \hfill $\Box$

\medskip

There is a direct connection between the total number of natural expansions and the total number of Fibonacci expansions.

\begin{theorem}\label{th:COUNT}
For every $N>3$ let $\beta(N) = d_{L(N)}\dots d_{R(N)}$ be the Bergman expansion of $N$. Then  $$\TotNu(N)=\TotFIB(F_{-R(N)+2}\,N).$$
\end{theorem}

\noindent  {\it Proof:} Suppose that $\beta(N) = d_{L}\dots d_R$, so $N=\sum_L^R d_i\varphi^i$.\\  Multiply by $\varphi^{-R+2}$:\\[-.1cm]
$$\varphi^{-R+2}N=\sum_{i=R}^L d_i\varphi^{i-R+2}= \sum_{j=2}^{L-R+2} d_{j+R-2}\varphi^j=  \sum_{j=2}^{L-R+2} e_j\varphi^j$$
where we substituted $j=i-R+2$, and defined $e_j:=j+R-2$.

Next we use the well known equation {$\varphi^j=F_j\varphi+F_{j-1}$}:
$$[F_{-R+2}\varphi+F_{-R+1}]N=\sum_{j=2}^{L-R+2} e_j[F_j\varphi+F_{j-1}].$$
This implies that\\[-.3cm]
$$F_{-R+2}N=\sum_{j=2}^{L-R+2} e_jF_j.$$
We conclude that  the number $F_{-R+2}N$ has  a Zeckendorf expansion given by the sum on the right side.

But the manipulations above can be made for any 0-1-word of length $L-R+1$, so the  golden mean flips of $d_L\dots d_R$  are in 1-to-1 correspondence with golden mean flips of $e_2\dots e_{L-R+2}$.
This implies that $\TotNu(N)=\TotFIB(F_{-R(N)+2}N).$ \hfill $\Box$

\bigskip

\noindent {\bf Example 1} The Bergman expansion of $4$ is $101\redcdot01$, and $F_4=3$. So $\TotNu(4)=\TotFIB(12)=1$.

\medskip

\noindent {\bf Example 2} The Bergman expansion of $14$ is $100100\redcdot001001$, and $F_8=21$.  So $\TotNu(14)=\TotFIB(294)=12$.

\medskip

\noindent {\bf Example 3}\label{ex:3}
Consider the Lucas numbers.  From $L_{2n}=\varphi^{2n}+\varphi^{-2n}$, and $L_{2n+1}=L_{2n}+L_{2n-1}$:\\
\hspace*{2cm} $\beta(L_{2n}) = 10^{2n}\redcdot0^{2n-1}1,\quad \beta(L_{2n+1}) = 1(01)^n\redcdot(01)^n.$\\
We read off:  $R(L_{2n})=-2n, R(L_{2n+1})=-2n$.

It is also clear that $\TotNu(L_{2n})=2n$, and $\TotNu(L_{2n+1})=1$.\\
So Theorem \ref{th:COUNT} gives the total number of Fibonacci representations of $F_{2n+2}L_{2n}$ and $F_{2n+2}L_{2n+1}$:

$\TotFIB(F_{2n+2}L_{2n}) = 2n, \;  \TotFIB(F_{2n+2}L_{2n+1}) = 1 $ for all $n\ge 1$.

\noindent We find in \cite{oeis}: From Miklos Kristof, Mar 19 2007:

  Let $L(n)$ = A000032($n$) = Lucas numbers. Then  for $a >= b$ and odd $b,\;  F(a+b) - F(a-b) = F(a)*L(b).$\\
 So $F_{2n+2}L_{2n+1}=F_{4n+3}-F_1=F_{4n+3}-1$. \; But $\TotFIB(F_n-1)=1$ is a well-known formula.

 \bigskip

\begin{remark}  An alternative proof of Theorem \ref{th:TotNuFn} can be given with Theorem \ref{th:COUNT}.

From Proposition \ref{prop:s} we  know that a number $N$  with $\beta(N)=d_L\ldots d_R$ in $\Lambda_{2n}$ has $-R(N)=2n$.
According to Lemma \ref{lem:FibLuc}: $F_{2n+2}\in \Lambda_{2n}$.  So Theorem \ref{th:COUNT}  leads to
    $$\TotNu(F_{2n+2})=\TotFIB(F_{2n+2}\,F_{2n+2})=\TotFIB(F_{2n+2}^2).$$
To finish the alternative proof, one needs to know that $\TotFIB(F_{2n}^2)=F_{2n-1}$ for all $n\ge 1$.
This can be proved similarly to the proof of the main result of Stockmeyers paper \cite{Stockmeyer}. The extra trick is to jointly prove the formula $\TotFIB(F_{2n}^2)=F_{2n-1}$ together with the formula $\TotFIB(F_{2n+1}^2-2)=F_{2n}.$

  The main tool of the proof is Klarner's result from \cite{Klarner66}:
   $\TotFIB(n)=\TotFIB(n-F_{m})+\TotFIB(F_{m+1}-n-2)$,\\
which holds for $m\ge 4$ and $F_m \le n < F_{m+1}-1$.

The proof by induction applies Klarner's identity with $n=F_m^2$, respectively   with $n=F_m^2-2$.\\
Here the identities (5): $F_m^2  =F_{2m-2}+F_{m-2}^2$ and (6): $F_m^2 =F_{2m-1}-F_{m-1}^2$ from \cite{Stockmeyer} make the induction step work.
\end{remark}

\section{Comparing Knott expansions and natural expansions}\label{sec:compare}

 It is not hard to see that the double golden mean shift---in general combined with more golden mean shifts---can be applied if and only if  the expansion ends in $10^s1$, where $s$ is odd. So the difference between the Knott expansions and the natural expansions is made more explicit by part a) of the following result.

\begin{proposition}\label{prop:s}{\bf a)}\: A number $N$ is in $\Lambda_{2n}$ if and only if $\beta(N)=...10^s1$, where $s$ is odd, and  $N$ is in $\Lambda_{2n+1}$ if and only if $\beta(N)=...10^s1$, where $s$ is even.\\
{\bf b)} \: Let $\beta(N)=L(N)...R(N)$.  A number $N$  in $\Lambda_{2n}$ has $-R(N)=2n$, a number $N$  in $\Lambda_{2n+1}$ has $-R(N)=2n+2$.
\end{proposition}

\medskip

Proposition \ref{prop:s} will be proved by induction. Thus we need recursions to let the proof work. These are given in the paper \cite{Dekk-how-to-add-FQ}, from which we repeat the following.

To obtain recursive relations, the interval $\Lambda_{2n+1}=[L_{2n+1}+1, L_{2n+2}-1]$ has to be divided into three subintervals. These three intervals are\\[-.4cm]
 \begin{align*}
I_n:=&[L_{2n+1}+1,\, L_{2n+1}+L_{2n-2}-1],\\
J_n:=&[L_{2n+1}+L_{2n-2},\, L_{2n+1}+L_{2n-1}],\\
K_n:=&[L_{2n+1}+L_{2n-1}+1,\, L_{2n+2}-1].
\end{align*}

\noindent It will be convenient to extend the monoid of  words of 0's and 1's to the corresponding free group. So, for example, $(01)^{-1}0001=1^{-1}001$.

\begin{theorem}{\bf [Recursive structure theorem, \cite{Dekk-how-to-add-FQ}]}\label{th:rec}

\noindent{\,\bf I\;} For all $n\ge 1$ and $k=0,\dots,L_{2n-1}$
one has $ \beta(L_{2n}+k) =  \beta(L_{2n})+ \beta(k) = 10\dots0 \,\beta(k)\, 0\dots 01.$\\
\noindent{\bf II\;} For all $n\ge 2$ and $k=1,\dots,L_{2n-2}-1$
\begin{align*}
I_n:&\quad \beta(L_{2n+1}+k) = 1000(10)^{-1}\beta(L_{2n-1}+k)(01)^{-1}1001,\\ K_n:&\quad\beta(L_{2n+1}+L_{2n-1}+k)=1010(10)^{-1}\beta(L_{2n-1}+k)(01)^{-1}0001.
\end{align*}
Moreover, for all $n\ge 2$ and $k=0,\dots,L_{2n-3}$
$$\hspace*{0.7cm}J_n:\quad\beta(L_{2n+1}+L_{2n-2}+k) = 10010(10)^{-1}\beta(L_{2n-2}+k)(01)^{-1}001001.$$
\end{theorem}

\medskip

\noindent {\it Proof of Proposition \ref{prop:s}:} To start the induction, we note that
\begin{align*}
  \Lambda_2 =[3,4];&\quad \beta(3)=100\redcdot01,\; \beta(4)=101\redcdot01, \\
   \Lambda_3 =[5,6];&\quad \beta(5)=1000\redcdot1001,\; \beta(6)=1010\redcdot0001.
\end{align*}

For the even intervals we have that $\beta(L_{2n}) = 10^{2n}\redcdot0^{2n-1}1$, so the expansion of the first element ends indeed in $10^s1$, where $s$ is odd. Note also that $R(L_{2n})=2n$, and this property will hold for all $L_{2n}+k$, $k=0,\dots,L_{2n-1}$ since the sum $\beta(L_{2n})+ \beta(k)$ in {\bf I} does not change the length of the negative part. Moreover, since  the length of the negative part of each $\beta(k)$ in the sum $\beta(L_{2n})+ \beta(k)$ is even (by the induction hypothesis for part {\bf b)}), the expansion must end in $10^s1$ with $s$ odd, simply because the difference of two even numbers is even.

\medskip

For the odd intervals we have to consider the three cases from {\bf II}.

For $I_n$: we know that $\beta(L_{2n-1}+k)$ ends in 01, so $\beta(L_{2n+1}+k)$ ends in 1001. For part {\bf b)}: the length of the negative part is increased by 2.

\smallskip

For $K_n$: $L_{2n-1}+k$ is from an odd interval, so the expansion ends in $10^{2t}1$ from some $t>0$. But then the expansion of $L_{2n+1}+L_{2n-1}+k$ ends in
$10^{2t}1\,(01)^{-1}0001=10^{2t-1}0001=10^{2t+2}1$. For part {\bf b)}: the length of the negative part is increased by 2.

\smallskip

For $J_n$: obviously $\beta(L_{2n+1}+L_{2n-2}+k)$ ends in 1001. For part {\bf b)}: the length of the negative part is $2n-2+ 4=2n+2$. \hfill $\Box$

\end{document}